\documentclass[12pt]{article}
\usepackage{amsfonts,amssymb,amsmath, epsfig}
\usepackage{color,graphicx,graphics,psfrag}
\usepackage{amsmath,amstext,amssymb,amsfonts, amscd}

\usepackage[numbers]{natbib}

\title{Crouzeix-Raviart MsFEM with Bubble Functions for Diffusion and Advection-Diffusion in Perforated Media}
\author{P. Degond$^{1,2}$, Alexei Lozinski$^{3}$, Bagus Putra Muljadi$^{1,2}$, Jacek Narski$^{1,2}$} 
\date{} 
\begin{document}
\maketitle

\vspace{0.5 cm}
\begin{center}
1-Universit\'e de Toulouse; UPS, INSA, UT1, UTM ;\\ 
Institut de Math\'ematiques de Toulouse ; \\
F-31062 Toulouse, France. \\
2-CNRS; Institut de Math\'ematiques de Toulouse UMR 5219 ;\\ 
F-31062 Toulouse, France.\\
email: pierre.degond@math.univ-toulouse.fr, muljadi.bp@gmail.com, jacek.narski@math.univ-toulouse.fr
\end{center}

\begin{center}
3-Laboratoire de Math\'ematiques, UMR CNRS 6623, Université de Franche-Comt\'e, 25030 Besançon Cedex, France\\
email: alexei.lozinski@univ-fcomte.fr
\end{center}

\begin{abstract}
The adaptation of Crouzeix - Raviart finite element in the context of multiscale finite element method (MsFEM) is studied and implemented on diffusion and advection-diffusion problems in perforated media. It is known that the approximation of boundary condition on coarse element edges when computing the multiscale basis functions critically influences the eventual accuracy of any MsFEM approaches. The weakly enforced continuity of Crouzeix - Raviart function space across element edges leads to a natural boundary condition for the multiscale basis functions which relaxes the sensitivity of our method to complex patterns of perforations. Another ingredient to our method is the application of bubble functions which is shown to be instrumental in maintaining high accuracy amid dense perforations. Additionally, the application of penalization method makes it possible to avoid complex unstructured domain and allows extensive use of simpler Cartesian meshes. 
\end{abstract}

\medskip
\noindent
{\bf Acknowledgements:} This work has been supported by the 'Fondation Sciences et Technologies pour l'Aéronautique et l'Espace', in the frame of the project 'AGREMEL' 
(contract \# RTRA-STAE/2011/AGREMEL/02). 


\section{Introduction}
Many important problems in modern engineering context have multiple-scale solutions e.g., transport in truly heterogeneous media like composite materials or in perforated media, or turbulence in high Reynolds number flows are some of the examples. Complete numerical analysis of these problems are difficult simply because they exhaust computational resources. In recent years, the world sees the advent of computational architectures such as parallel and GPU programming; both are shown to be advantageous to tackle resource demanding problems. Nevertheless, the size of the discrete problems remains big. In some engineering contexts, it is sometimes sufficient to predict macroscopic properties of multiscale systems. Hence it is desirable to develop an efficient computational algorithm to solve multiscale problems without being confined to solving fine scale solutions. Several methods sprung from this purpose namely, Generalized finite element methods \citep{babuskaetal1}, wavelet-based numerical homogenization method \citep{dorobantuengquist}, variational multiscale method \citep{nolenetal},various methods derived from homogenization theory \citep{bourgeat}, equation-free computations \citep{kevrekidisetal}, heterogeneous multiscale method \citep{weinanengquist} and many others. In the context of diffusion in perforated media, some studies have been done both theoretically and numerically in \citep{cioranescuetal},\citep{cioranescumurat},\citep{henningohlberger},\citep{hornung}, and \citep{lions}. For the case of advection-diffusion a method derived from heterogeneous multiscale method addressing oscillatory coefficients is studied in \citep{Deng20091549}.

In this paper, we present the development of a dedicated solver for solving multiscale problems in perforated media most efficiently. We confine ourselves in dealing with only stationary diffusion and advection-diffusion problems as means to pave the way toward solving more complicated problems like Stokes. We begin by adapting the concept of multiscale finite element method (MsFEM) originally reported in \citep{thomhou}. The MsFEM method relies on the expansion of the solution on special basis functions which are pre-calculated by means of local simulations on a fine mesh and which model the microstructure of the problem. By contrast to sub-grid modeling methodologies, the multiscale basis functions are calculated from the actual geometry of the domain and do not depend on an (often arbitrary) analytical model of the microstructure. A study on the application of MsFEM in porous media has been done in \citep{Efendiev2007577}, and although it could have bold significance in geo- or biosciences, they can be applied also in different contexts, e.g., pollutant dispersion in urban area \citep{laetitia} or on similar problems which are extremely dependent on the geometry of perforations but their full account leads to very time consuming simulations. Textbook materials on the basics of MsFEM can be found in \citep{houefendiev}.

It is understood that when constructing the multiscale basis function, the treatments of boundary condition on coarse elements greatly influence the accuracy of the method of interest. For example, in the original work of Hou and Wu, the oversampling method was introduced to provide the best approximation of the boundary condition of the multiscale basis functions which is also of high importance when dealing with non-periodic perforations.  Oversampling here means that the local problem in the course element are solved on a domain larger than the element itself, but only the interior information is communicated to the coarse scale equation. This reduces the effect of wrong boundary conditions and bad sampling sizes. The ways in which the sampled domain is extended lead to various oversampling methods, see \citep{houefendiev}, \citep{chuetal}, \citep{henningpeterseim}, \citep{efendievetal}. The non-conforming nature of Crouzeix-Raviart element, see \citep{CRRairo}, is shown to provide great 'flexibility' especially when non-periodically perforated media is considered. In the construction of Crouzeix-Raviart multiscale basis functions, the conformity between coarse elements are not enforced in a strong sense, but rather in a weak sense i.e., the method requires merely the average of the ''jump'' of the function to vanish at coarse element edges. When very dense perforations are introduced, which often makes it virtually impossible to avoid intersections between coarse element edges and perforations, the benefit of using Crouzeix-Raviart MsFEM is significant for it allows the multiscale basis functions to have natural boundary conditions on element edges making it insensitive to complex patterns of perforations. Moreover, the integrated application of penalization method enables one to carry the simulations onto simple Cartesian meshes. Note that some methods derived from homogenization theory may provide robust and accurate results provided that the underlying multiscale structure or subgrid effects satisfies the necessary constraints which is not the case for problems with non-periodic perforations.  In this paper several computational results with non-periodic perforations are given to highlight the feasibility of our method in such circumstances. Another important ingredient to our method is the multiscale finite element space enrichment with bubble functions. Again, when very dense perforations are considered, it is both crucial and difficult to capture correct approximations between perforations for which the application of bubble functions is offered as the remedy. We illustrate these problems in our paper to highlight the contribution of bubble function in improving the accuracy of our MsFEM.  Our work continues the application of Crouzeix-Raviart MsFEM done on prototypical elliptic problems \citep{clbetal} and on diffusion problems with homogeneous boundary condition \citep{lozinskibubbleetal}. Improvements are done to the earlier work by introducing \textit{bubble functions} and to the latter by extending the application to advection-diffusion problems with non-homogeneous boundary conditions. 

The paper is organized as the following. In chapter \ref{sec:crouzeix} we begun with the formulation of the problem and the construction of our MsFEM. Here we explain the construction of Crouzeix-Raviart MsFEM functions space with bubble functions and the multiscale basis functions. In chapter \ref{sec:boundary} the application of non-homogeneous boundary conditions is explained. In chapter \ref{sec:applicationpenal} we  describe the application of penalization method. Demonstrations of our MsFEM in terms of computational simulations and its analysis can be found in chapter \ref{sec:numerical} followed by some concluding remarks.

\section{Crouzeix-Raviart MsFEM with bubble functions enrichment}
\label{sec:crouzeix}
We consider an advection-diffusion problem laid in a bounded domain $\Omega \in \mathcal{R}^d$ within which a set $B_\epsilon$ of perforations is included. From here on we assume that the ambient dimension is $d = 2$. The perforated domain with voids left by perforations is denoted $\Omega^\epsilon = \Omega \setminus B_\epsilon$, where $\epsilon$ denotes the minimum width of perforations. The advection-diffusion problem is then to find $u : \Omega^\epsilon \rightarrow \mathcal{R}$ which is the solution to
\begin{eqnarray}
 -\nabla \cdot (\mathcal{A} \nabla u) + \vec{w} \cdot \nabla u = f & \textrm{in} & \Omega^\epsilon \label{maineq}\\
 u = 0 & \textrm{on} & \partial B^\epsilon \cap \partial\Omega^\epsilon \nonumber \\
 u = g & \textrm{on} & \partial \Omega \cap \partial\Omega^\epsilon \nonumber
\end{eqnarray}
where $f : \Omega \rightarrow \mathcal{R}$ is a given function, $g$ is a function fixed on boundary $\partial\Omega$ and $\vec{w}$ is a given velocity field. In this paper, we consider only the Dirichlet boundary condition on $\partial B_\epsilon$ namely $u_{|\partial B^\epsilon} = 0$ thereby assuming that the perforation is opaque. Other kinds of boundary conditions on $\partial B_\epsilon$ are subject to a completely new endeavour. Recent works on Crouzeix-Raviart MsFEM focusing on diffusion problems with homogeneous boundary condition $g = 0$ has been done in \citep{lozinskibubbleetal}.

When linear boundary condition for MsFEM basis function is considered, it is difficult to approximate correct coarse node solutions when one or more perforations coincide with any of the coarse element's boundaries. The approximation of the MsFEM basis function will be distorted and the whole solution will be affected. This problem often and can be relaxed by using oversampling methods. However, in practice, this brings an inherent inconvenience since the size and position of perforations are most of the times unpredictable which requires some problem dependent parameters to be introduced. Moreover, the computational cost also increases due to enlarged sampled domain. 

The Crouzeix-Raviart basis functions are non-conforming throughout the computational domain. The continuity of the functions are enforced weakly i.e., it requires no fixed values across the boundaries but rather vanishing ''jump'' averages on each edge. In order to explain the MsFEM space in the vein of Crouzeix-Raviart's finite element, we define a mesh $\mathcal{T}_H$ in $\Omega$ which are discrete polygons with each diameter at most $H$ and made up by $n_H$ coarse elements and $n_e$ coarse element edges. Denote $\epsilon_H$ the set of all edges $e$ of $\mathcal{T}_H$ which includes edges on the domain boundary $\partial\Omega$. It is assumed that the mesh does not include any hanging nodes and each edge is shared by two elements except those on $\partial\Omega$ which belongs only to one element. $\mathcal{T}_H$ is assumed a regular mesh. By regular mesh, we mean for any mesh element $T \in \mathcal{T}_H$, there exists a smooth one-to-one mapping $\mathcal{M} : \bar{T} \rightarrow T$ where $\bar{T} \subset \mathcal{R}^d$ is the element of reference, and that $\parallel \nabla\mathcal{M} \parallel_{L^\infty} \leq DH$, $\parallel \nabla\mathcal{M}^{-1} \parallel_{L^\infty} \leq DH^{-1}$ with $D$ being universal constant independent of $T$. We introduce the functional space for Crouzeix-Raviart type MsFEM with bubble function enrichment 
\begin{eqnarray}
V_H &=& \{u \in L^2(\Omega) \hspace*{1mm} | \hspace*{1mm} u_{|T}\in H^1(T) \hspace*{1mm}\textrm{for all}\hspace*{1mm} T\in\mathcal{T}_H, \label{crspace}\nonumber\\
&&-\nabla \cdot (\mathcal{A}\nabla u) + \vec{w}\cdot\nabla u= \textrm{constant in}\hspace*{1mm} T \cap \Omega^\epsilon \hspace*{1mm}\textrm{for all}\hspace*{1mm} T\in\mathcal{T}_H,  \nonumber\\ 
&&u = 0 \hspace*{1mm}\textrm{on}\hspace*{1mm} \partial B_\epsilon,\hspace*{1mm}n\cdot\nabla u = \hspace*{1mm}\textrm{constant on}\hspace*{1mm} e \cap \Omega^\epsilon \hspace*{1mm}\textrm{for all}\hspace*{1mm} e\in\epsilon_H,\nonumber\\
&&\int_{e}[[u]]=0\hspace*{1mm}\textrm{for all}\hspace*{1mm} e \in \epsilon_H \footnotemark\}
\end{eqnarray}
\footnotetext{ $ \textrm{On}\hspace*{1mm} e\in\epsilon_H\cap\partial\Omega, \hspace*{1mm}\hspace*{1mm}\int_{e}[[u]] = \int_{e}u. $  }
where $[[u]]$ denotes the jump of $u$ over an edge. 
The MsFEM approximation to Eq. (\ref{maineq}) is the solution of $u_H \in V_H$ to
\begin{equation}
a_H(u_H,v_H) = \int_{\Omega^\epsilon} f v_H \hspace{10mm} \textrm{for all}\hspace{10mm} v_H\in V_H
\end{equation}
where 
\begin{equation}
a_H(u,v) = \sum_{T\in\mathcal{T}_H} \left(\int_{T\cap\Omega^\epsilon}\nabla u \cdot \nabla v + \int_{T\cap\Omega^\epsilon} (\vec{w}\cdot\nabla u) v\right).
\end{equation}
The basis for $V_H$ contains functions associated to edges $e$ and mesh elements $T \cap (\Omega \setminus B_\epsilon) \cap \Omega^\epsilon$. The former has the notation $\Phi_e$ and the latter $\Phi_B$. The edges composing $T^k$ are denoted $\Gamma_i^k$ with $i = 1,\cdots,N_\Gamma$ whereas $k = 1, \cdots, n_H$. The Crouzeix-Raviart multiscale basis functions ${\Phi_e}^k_i$ are then the unique solution in $H^1(T^k)$ to
\begin{eqnarray}
-\nabla \cdot [\mathcal{A} \nabla {\Phi_e}^k_i] + \vec{w}\cdot\nabla{\Phi_e}^k_i &= 0& \, \textrm{in} \hspace{1cm}  T^k \\
\int_{\Gamma_i^k} {\Phi_e}^k_i &= \delta_{ie} &\textrm{for}\hspace{1cm}  i = 1,\cdots, N_\Gamma\\
n\cdot \mathcal{A} \nabla {\Phi_e}^k_i &=\lambda_i^k& \, \textrm{on}\hspace{1cm} \Gamma_i^k, i = 1,\cdots , N_\Gamma.
\end{eqnarray}
Whereas the bubble functions ${\Phi_B}^k$ can be obtained by solving for each element $T^k$
\begin{eqnarray}
-\nabla \cdot (\mathcal{A}\nabla {\Phi_B}^k)+ \vec{w}\cdot\nabla{\Phi_B}^k &= 1,& \hspace{1cm}\textrm{in}\hspace{1cm} T^k \nonumber\\
{\Phi_B}^k &= 0,& \hspace{1cm}\textrm{on}\hspace{1cm} \partial T^k
\end{eqnarray}
such that the approximated solution $u_H$ is described as
\begin{equation}
\label{disc.u}
u_H (x,y) = \sum_{i = 1}^{n_H}u_i {\Phi_e}_i (x,y) + \sum_{k = 1}^{n_e}u^k {\Phi_B}^k (x,y).
\end{equation}
%%%%%%%%%%%%%%%%%%%%%%%%%%%%%%%%%%%%%%%%%%%%%%%%%%%%%%%%%%%%%%%%%%%%%%%%%%%%%%%%%%
%%%%%%%%%%%%%%%%%%%%%%%%%%%%%%%%%%%FIGURES%%%%%%%%%%%%%%%%%%%%%%%%%%%%%%%%%%%%%%%%
%%%%%%%%%%%%%%%%%%%%%%%%%%%%%%%%%%%%%%%%%%%%%%%%%%%%%%%%%%%%%%%%%%%%%%%%%%%%%%%%%%
\begin{figure}[htbp]
\centering
\includegraphics[bb=0 0 3004 2900,width = 7in]{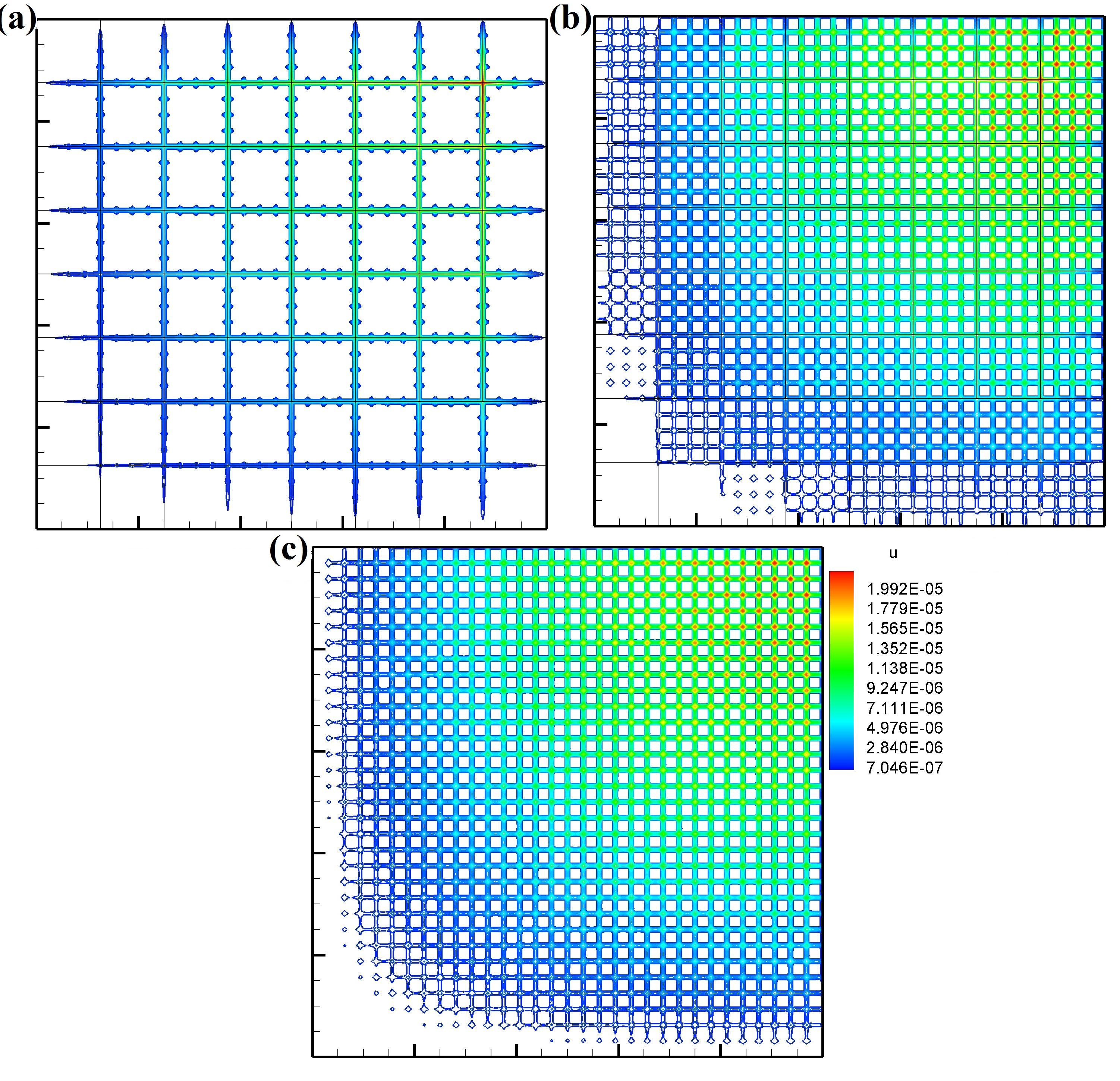} 
\caption{(a) $8 \times 8$ coarse elements without bubble functions, (b) $8 \times 8$ coarse elements with bubble functions compared with (c) Reference solution with $1024 \times 1024$ elements.}
\label{bubble}
\end{figure} 
%%%%%%%%%%%%%%%%%%%%%%%%%%%%%%%%%%%%%%%%%%%%%%%%%%%%%%%%%%%%%%%%%%%%%%%%%%%%%%%%%%
\begin{figure}[htbp]
\centering
\includegraphics[bb=0 0 1248 542,width = 7in]{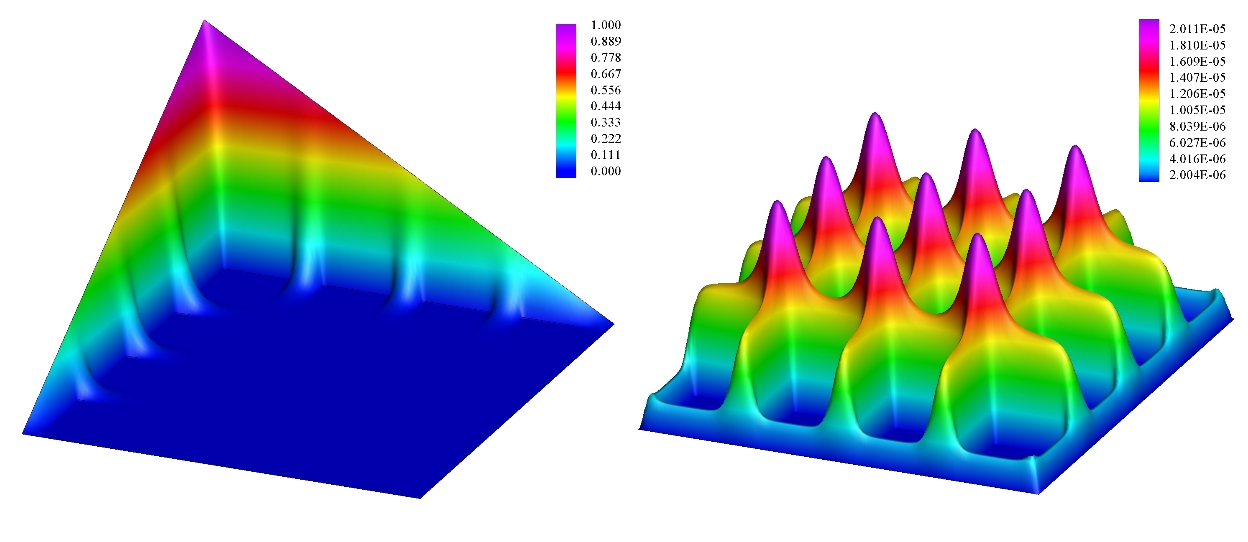} 
\caption{(a) Standard MsFEM basis function with (b) A bubble function in a coarse element.}
\label{bubblebf}
\end{figure} 
%%%%%%%%%%%%%%%%%%%%%%%%%%%%%%%%%%%%%%%%%%%%%%%%%%%%%%%%%%%%%%%%%%%%%%%%%%%%%%%%%%
\begin{figure}[htbp]
\centering
\includegraphics[bb=0 0 766 677,width = 5in]{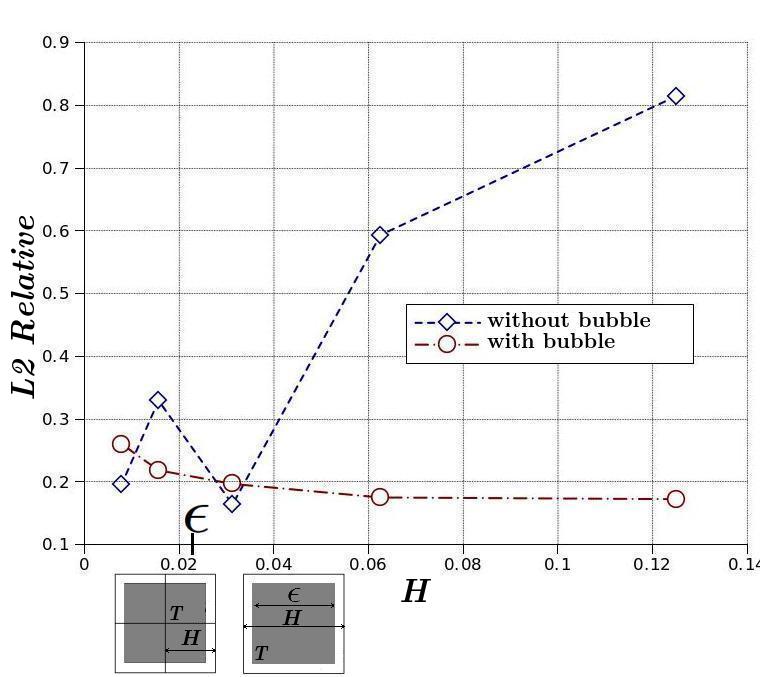} 
\caption{Relative error of standard MsFEM with and without bubble functions with respect to reference solution, $\epsilon = 0.021875$.}
\label{performance}
\end{figure} 
%%%%%%%%%%%%%%%%%%%%%%%%%%%%%%%%%%%%%%%%%%%%%%%%%%%%%%%%%%%%%%%%%%%%%%%%%%%%%%%%%%
\begin{figure}[htbp]
\centering
\includegraphics[bb=0 0 3181 2931,width = 7in]{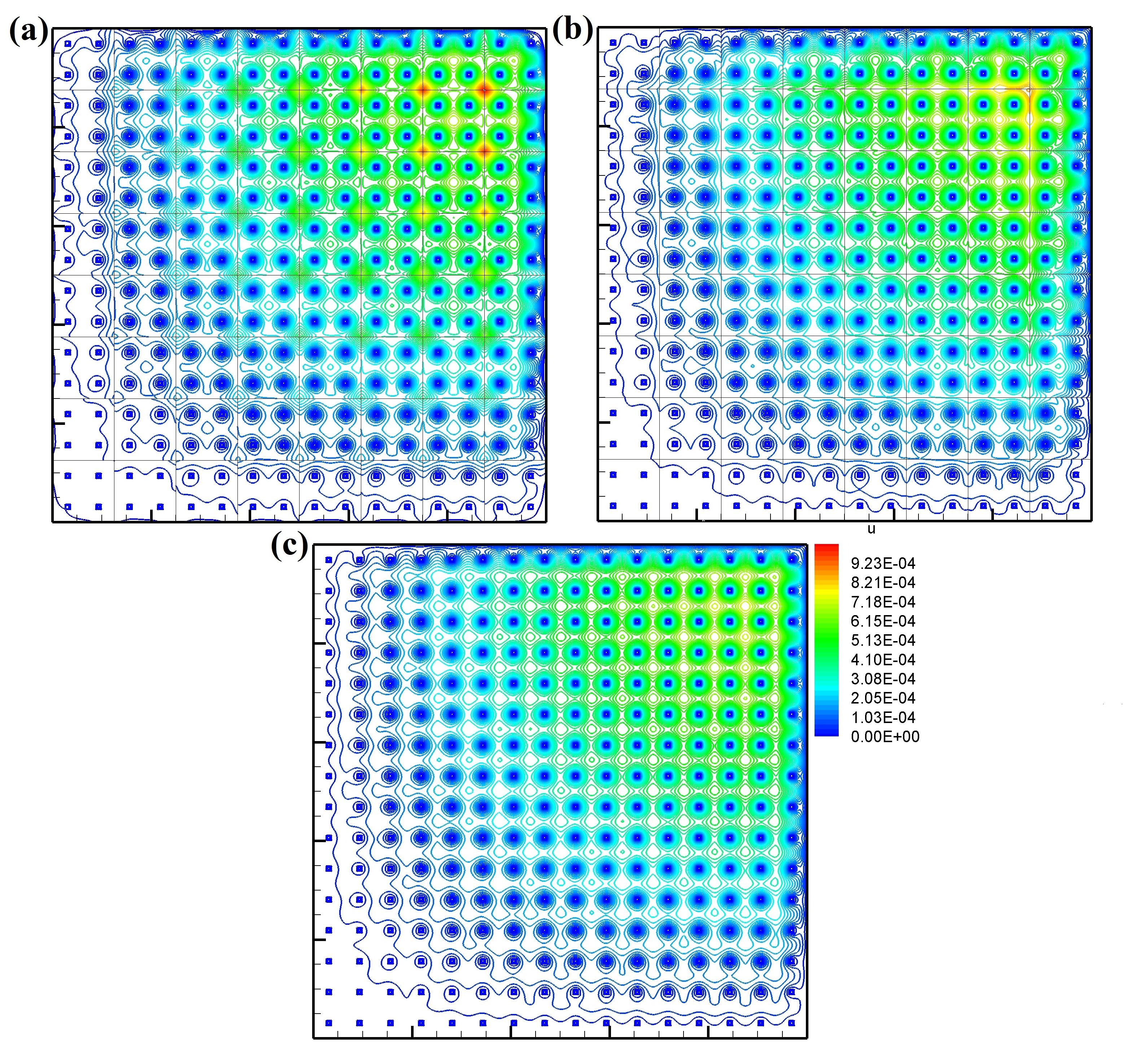} 
\caption{Coarse mesh-perforations-nonintersecting case solved on $8 \times 8$ coarse elements with: (a) Crouzeix-Raviart MsFEM with bubble functions (b) Standard MsFEM with bubble functions compared with (c) Q1 FEM solution as reference with $1024 \times 1024$ elements.}
\label{nonintersecting}
\end{figure} 
%%%%%%%%%%%%%%%%%%%%%%%%%%%%%%%%%%%%%%%%%%%%%%%%%%%%%%%%%%%%%%%%%%%%%%%%%%%%%%%%%%
\begin{figure}[htbp]
\centering
\includegraphics[bb=0 0 3271 2901,width = 7in]{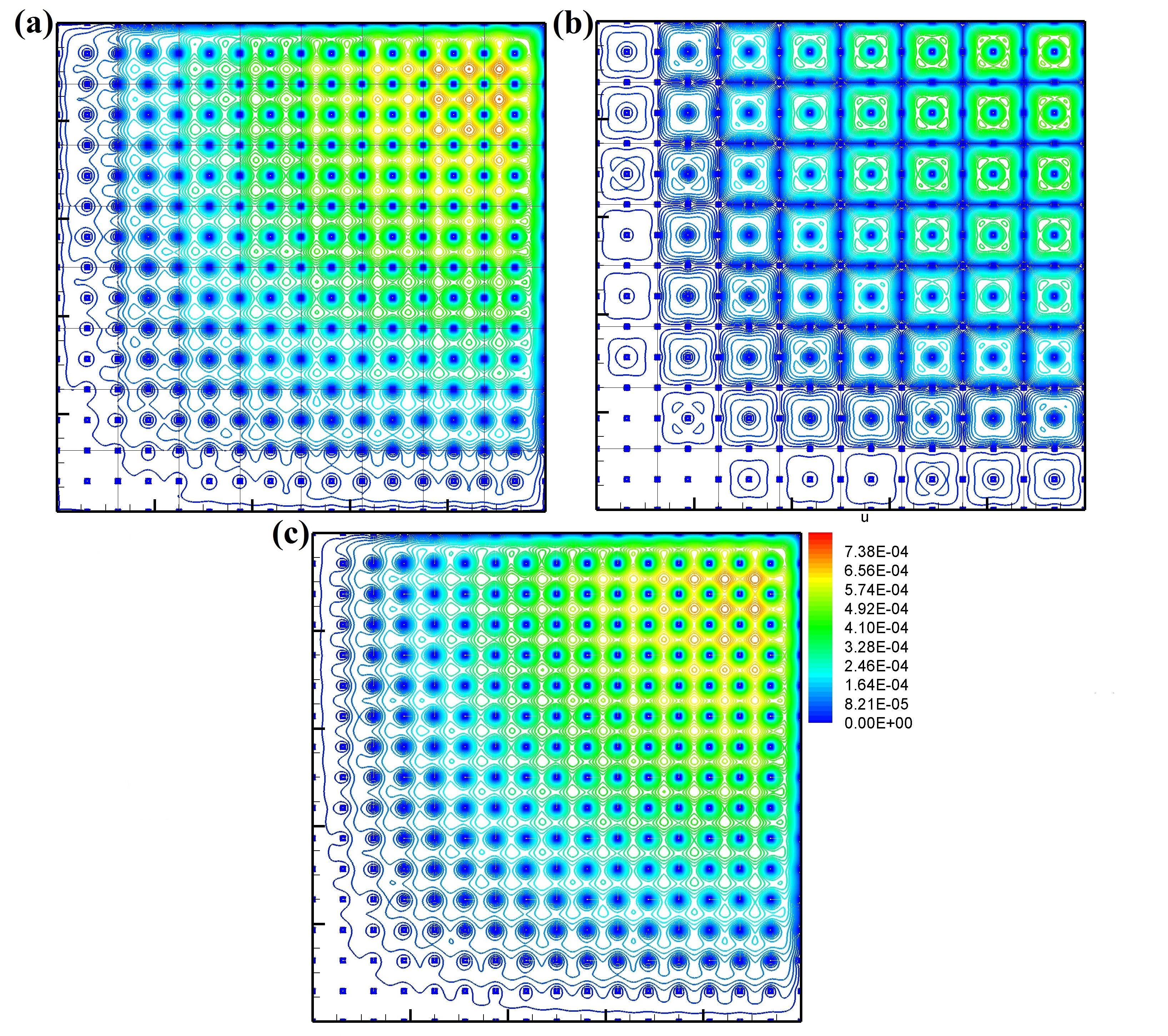} 
\caption{Coarse mesh-perforations-intersecting case on $8 \times 8$ coarse elements solved with: (a) Crouzeix-Raviart MsFEM with bubble functions (b) Standard MsFEM with bubble functions compared with (c) Q1 FEM solution as reference with $1024 \times 1024$ elements.}
\label{intersecting}
\end{figure} 
%%%%%%%%%%%%%%%%%%%%%%%%%%%%%%%%%%%%%%%%%%%%%%%%%%%%%%%%%%%%%%%%%%%%%%%%%%%%%%%%%%
%%%%%%%%%%%%%%%%%%%%%%%%%%%%%%%%%%%%%%%%%%%%%%%%%%%%%%%%%%%%%%%%%%%%%%%%%%%%%%%%%%
\begin{figure}[htbp]
\centering
\includegraphics[bb=0 0 3056 1688,width = 7in]{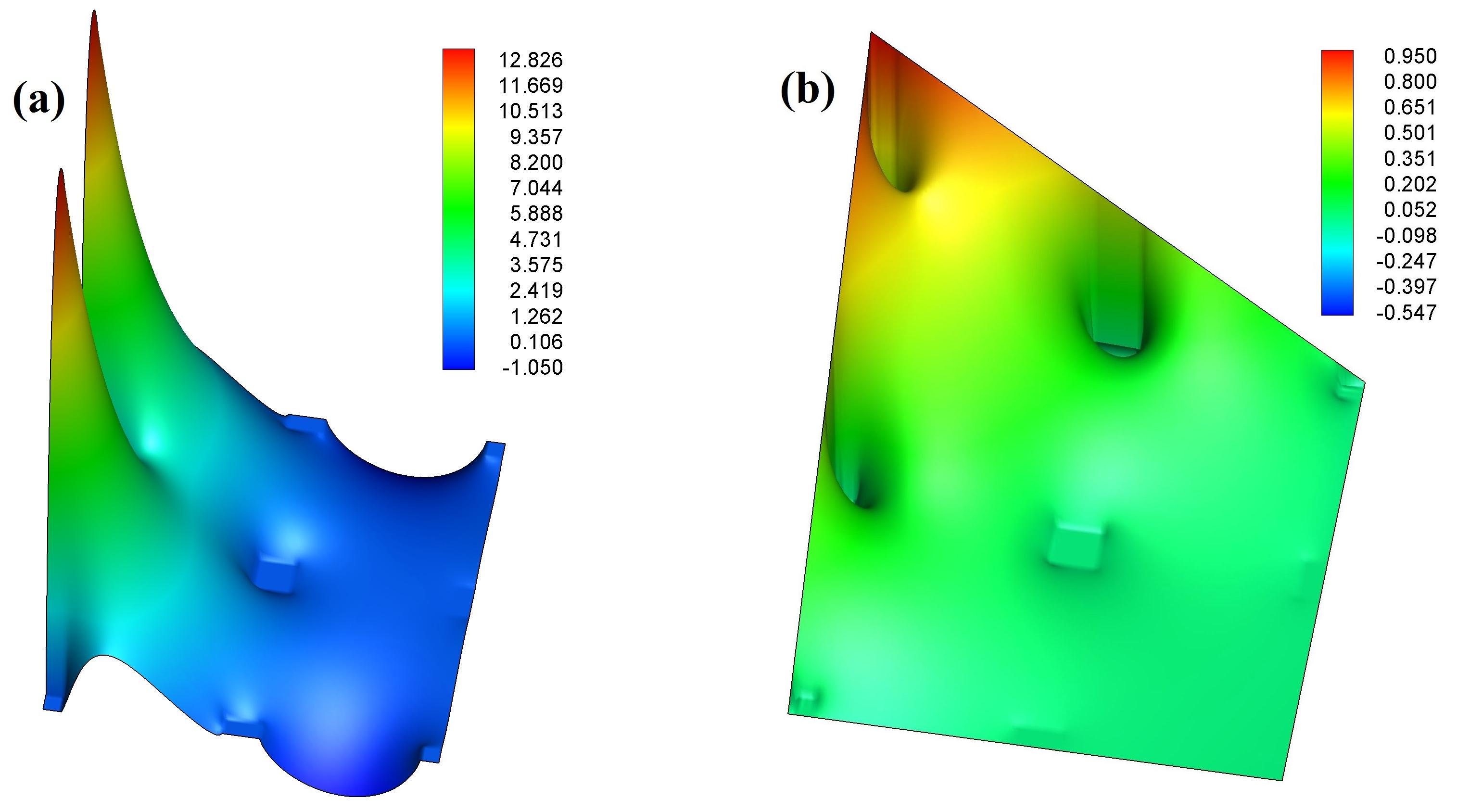} 
\caption{(a) Crouzeix-Raviart MsFEM basis function, (b) Nodal based MsFEM basis function without oversampling.}
\label{crbf}
\end{figure}
%%%%%%%%%%%%%%%%%%%%%%%%%%%%%%%%%%%%%%%%%%%%%%%%%%%%%%%%%%%%%%%%%%%%%%%%%%%%%%%%%%
\begin{figure}[htbp]
\centering
\includegraphics[bb = 0 0 3151 2921,width = 7in]{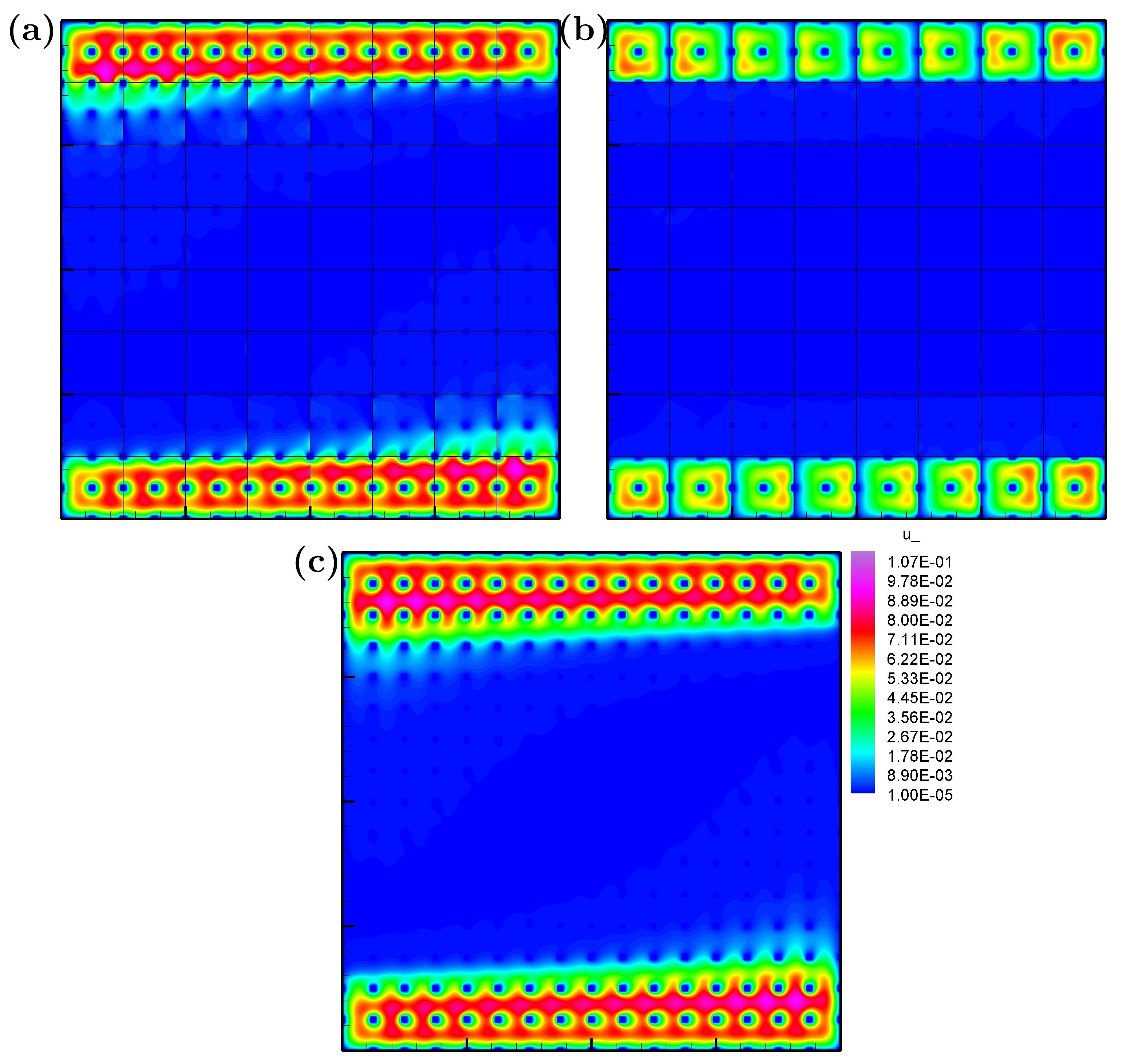} 
\caption{Coarse mesh-perforations-intersecting advection-diffusion case on $8 \times 8$ coarse elements solved with: (a) Crouzeix-Raviart MsFEM with bubble functions (b) Standard MsFEM with bubble functions and (c) Q1 FEM solution as reference with $1024 \times 1024$ elements, all with $\mathcal{A} = 0.03$.}
\label{convdiff}
\end{figure} 
%%%%%%%%%%%%%%%%%%%%%%%%%%%%%%%%%%%%%%%%%%%%%%%%%%%%%%%%%%%%%%%%%%%%%%%%%%%%%%%%%%
%%%%%%%%%%%%%%%%%%%%%%%%%%%%%%%%%%%%%%%%%%%%%%%%%%%%%%%%%%%%%%%%%%%%%%%%%%%%%%%%%%
\begin{figure}[htbp]
\centering
\includegraphics[bb = 0 0 1005 388,height= 3in]{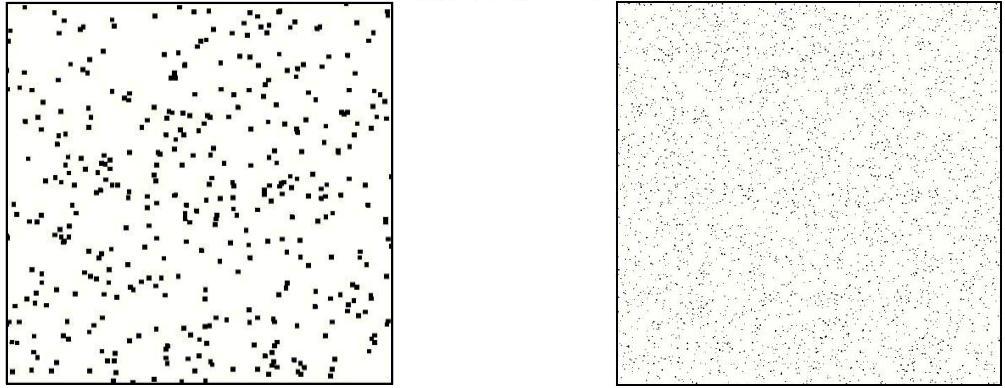} 
\caption{Non-periodic perforations: (left) Case a, (right) Case b}
\label{perfor}
\end{figure} 
%%%%%%%%%%%%%%%%%%%%%%%%%%%%%%%%%%%%%%%%%%%%%%%%%%%%%%%%%%%%%%%%%%%%%%%%%%%%%%%%%%
%%%%%%%%%%%%%%%%%%%%%%%%%%%%%%%%%%%%%%%%%%%%%%%%%%%%%%%%%%%%%%%%%%%%%%%%%%%%%%%%%
\begin{figure}[htbp]
\centering
\includegraphics[bb = 0 0 3141 2481,height = 7in]{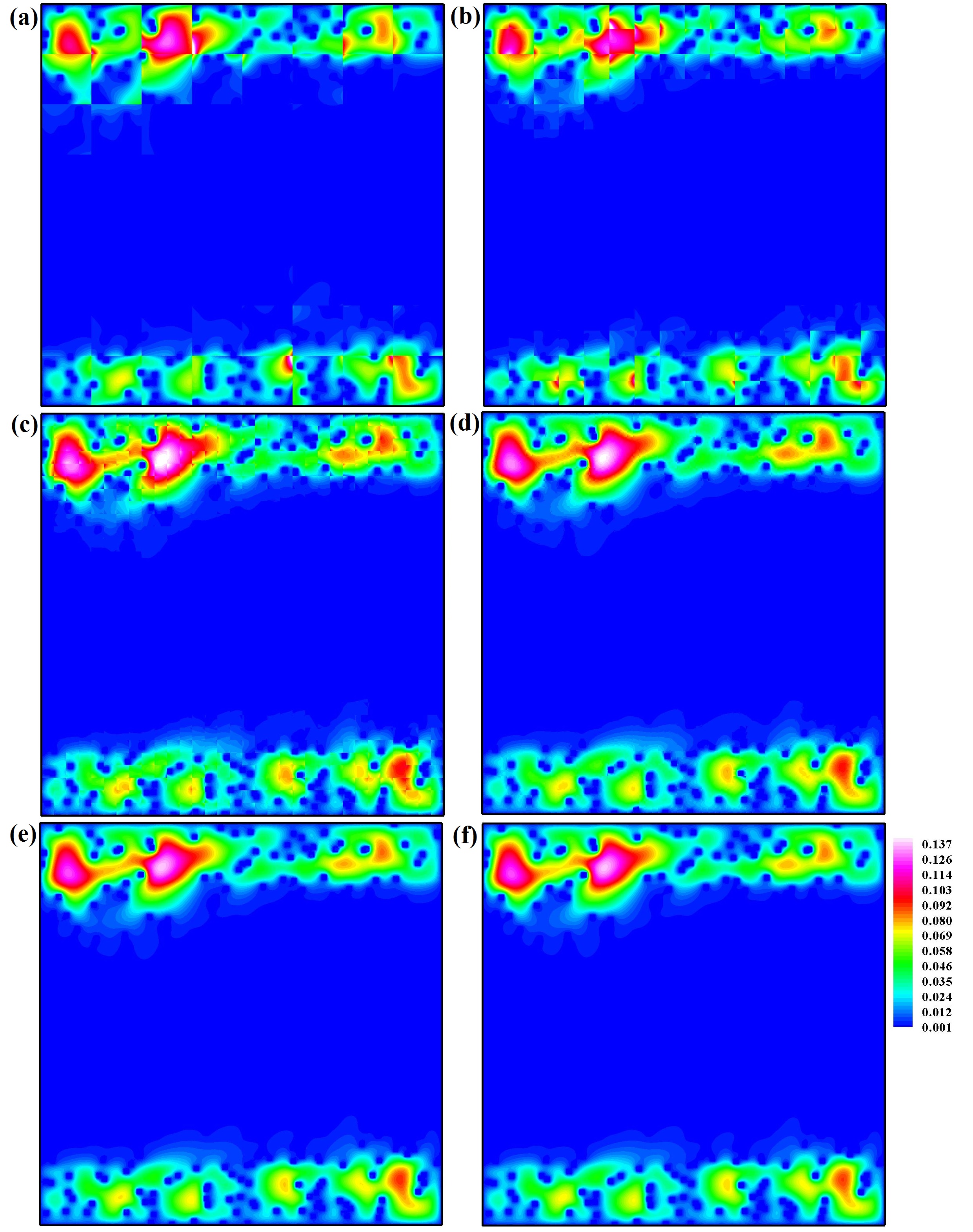} 
\caption{Advection-diffusion on non-periodically perforated domain (case a) (a)$8\times 8$, (b)$16\times 16$, (c)$32\times 32$, (d)$64\times 64$, (e)$128\times 128$, (f) Reference solution, Q1-Q1 FEM on 1024 $\times$ 1024 elements}
\label{crnp1}
\end{figure} 
%%%%%%%%%%%%%%%%%%%%%%%%%%%%%%%%%%%%%%%%%%%%%%%%%%%%%%%%%%%%%%%%%%%%%%%%%%%%%%%%%%
%%%%%%%%%%%%%%%%%%%%%%%%%%%%%%%%%%%%%%%%%%%%%%%%%%%%%%%%%%%%%%%%%%%%%%%%%%%%%%%%%%
\begin{figure}[htbp]
\centering
\includegraphics[bb = 0 0 3141 2481,height = 7in]{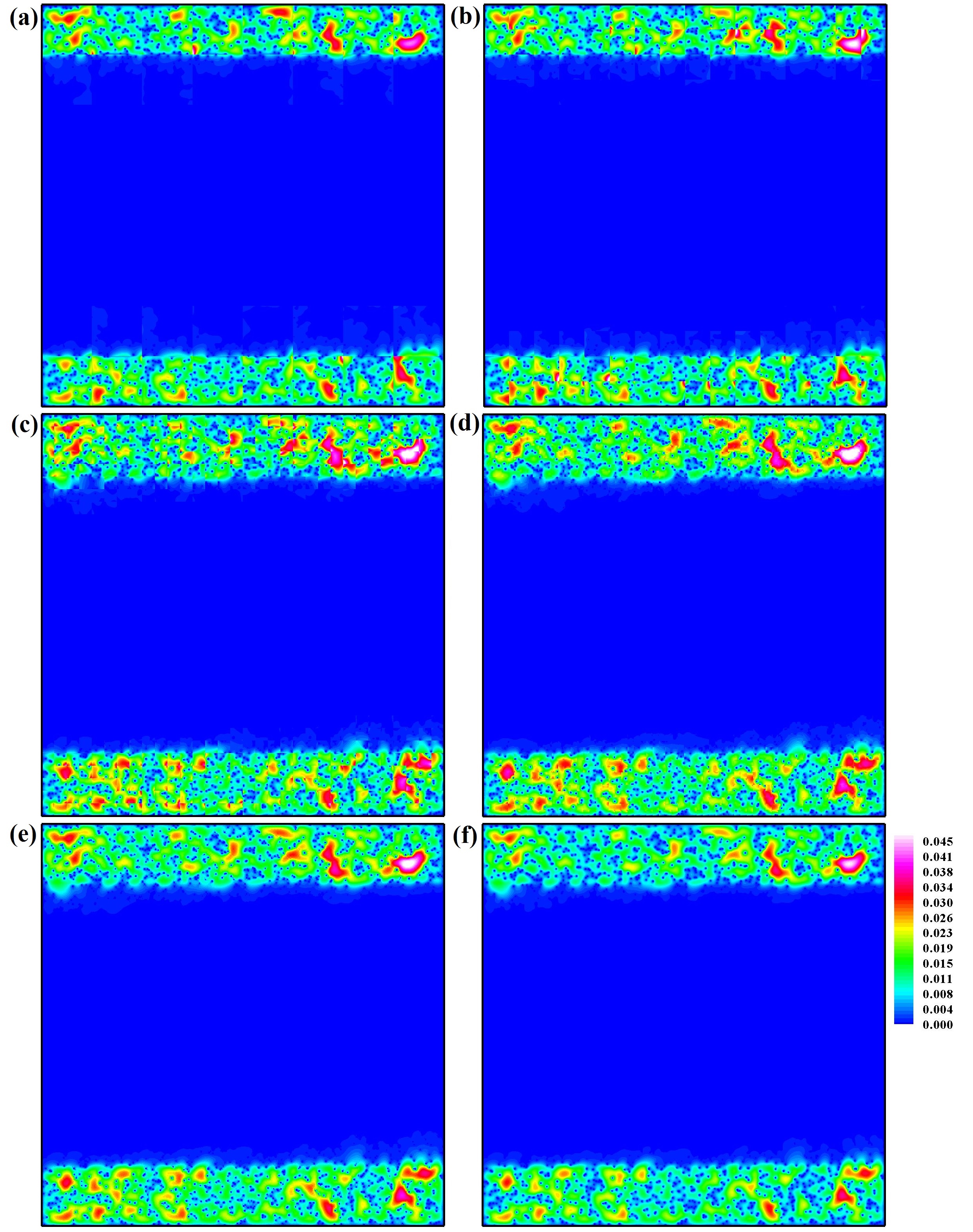} 
\caption{Advection-diffusion on non-periodically perforated domain (case b) (a)$8\times 8$, (b)$16\times 16$, (c)$32\times 32$, (d)$64\times 64$, (e)$128\times 128$, (f) Reference solution, Q1-Q1 FEM on 1024 $\times$ 1024 elements}
\label{crnp2}
\end{figure} 
%%%%%%%%%%%%%%%%%%%%%%%%%%%%%%%%%%%%%%%%%%%%%%%%%%%%%%%%%%%%%%%%%%%%%%%%%%%%%%%%%%
%%%%%%%%%%%%%%%%%%%%%%%%%%%%%%%%%%%%%%%%%%%%%%%%%%%%%%%%%%%%%%%%%%%%%%%%%%%%%%%%%%
\begin{figure}[htbp]
\centering
\includegraphics[bb = 0 0 379 382,height = 3in]{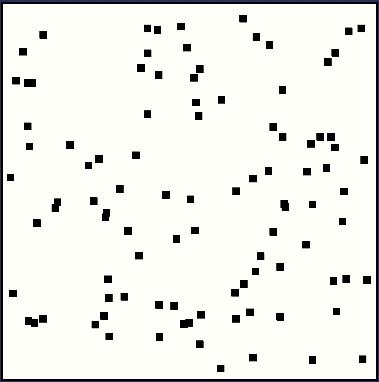} 
\caption{Domain with non-periodic perforations}
\label{perfor2}
\end{figure} 
%%%%%%%%%%%%%%%%%%%%%%%%%%%%%%%%%%%%%%%%%%%%%%%%%%%%%%%%%%%%%%%%%%%%%%%%%%%%%%%%%%
%%%%%%%%%%%%%%%%%%%%%%%%%%%%%%%%%%%%%%%%%%%%%%%%%%%%%%%%%%%%%%%%%%%%%%%%%%%%%%%%%%
\begin{figure}[htbp]
\centering
\includegraphics[bb = 0 0 3141 2481,height = 7in]{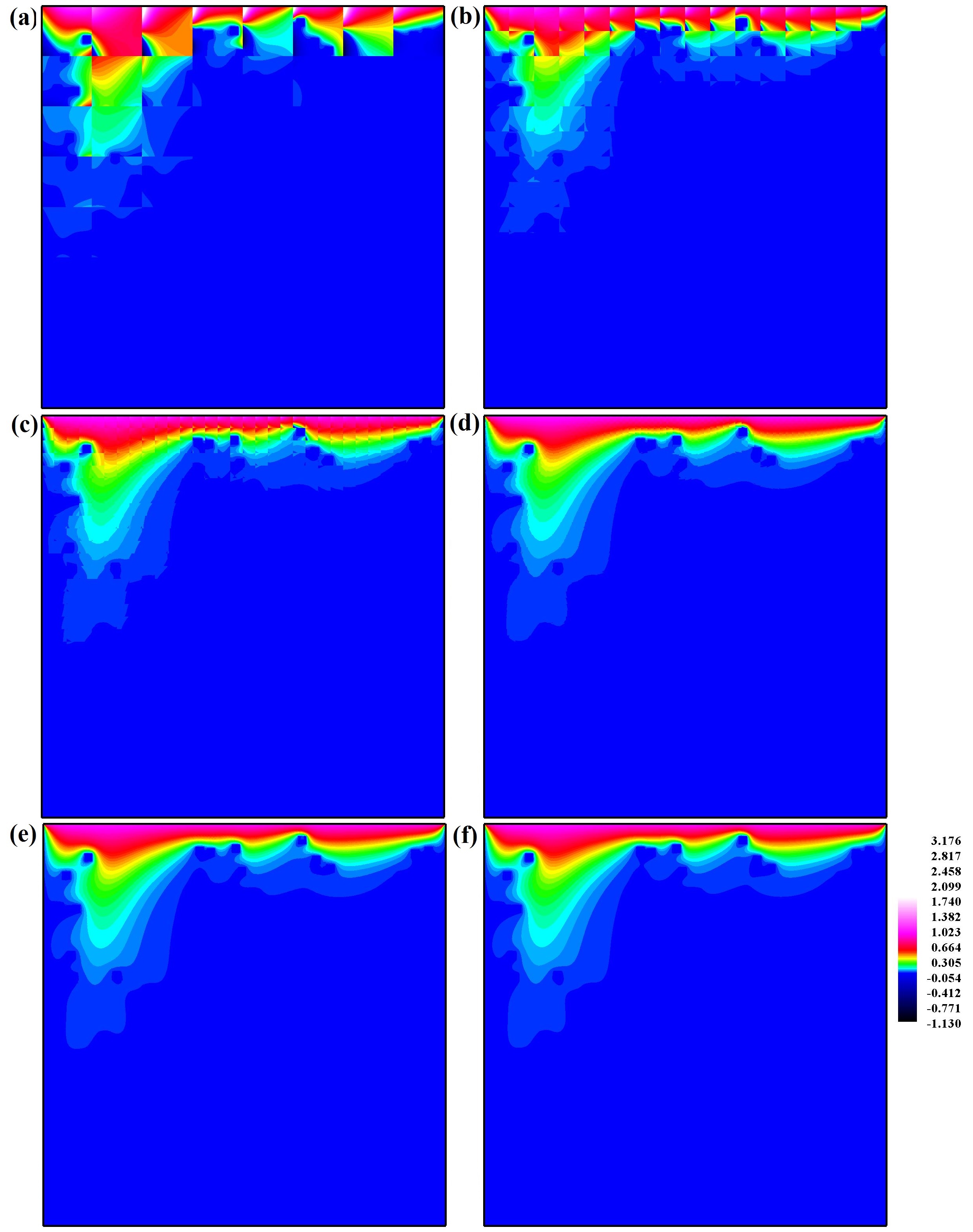} 
\caption{Advection-diffusion with non-periodic perforations and non-homogeneous boundary condition (a)$8\times 8$, (b)$16\times 16$, (c)$32\times 32$, (d)$64\times 64$, (e)$128\times 128$, (f) Reference solution, Q1-Q1 FEM on 1024 $\times$ 1024 elements}
\label{crnp3}
\end{figure} 
%%%%%%%%%%%%%%%%%%%%%%%%%%%%%%%%%%%%%%%%%%%%%%%%%%%%%%%%%%%%%%%%%%%%%%%%%%%%%%%%%%
In this paper, although the general formulation is focused on advection-diffusion problem, various tests on diffusion-only cases are presented for showing the contributions of bubble functions and Crouzeix-Raviart multiscale basis functions. The MsFEM formulation for diffusion cases is largely similar to the advection-diffusion counterpart. Regarding the advection-diffusion problems, one could consider another approach for selecting the space for the basis functions namely by following the Petrov-Galerkin formulation which could be useful for problems with very high P\'eclet numbers, nevertheless this is not the emphasis of this paper. Rigorous studies on Crouzeix-Raviart MsFEM's numerical analysis and error estimates for highly oscillatory elliptic problems and for diffusion problems in perforated media can be found in \citep{clbetal} and \citep{lozinskibubbleetal}

\section{Boundary condition}
\label{sec:boundary}
We propose to approximate the non-homogeneous dirichlet boundary condition in Eq. (\ref{maineq}) by,
\begin{equation}
\label{specialbc}
\int_{e \cap \partial\Omega} {u_H} = \int_{e\cap \partial\Omega} {g},\hspace{5mm}\textrm{for}\hspace{1mm}\textrm{all}\hspace{1mm} e \in \epsilon_H \hspace{1mm}\textrm{on}\hspace{1mm}\partial\Omega
\end{equation}
Equation (\ref{specialbc}) is therefore equivalent with 
\begin{equation}
u_{e\cap \partial\Omega} = \frac{1}{|e|}\int_{e\cap \partial\Omega} {g}.
\end{equation} 
With this, the construction of Crouzeix-Raviart basis functions associated both on edges at domain boundary or within the domain can be carried out in a similar fashion. This approach is a modification with respect to the earlier works in \citep{clbetal},\citep{lozinskibubbleetal} where the boundary condition were strongly incorporated in the definition of $V_H$.  Our approach therefore gives more flexibility when implementing non zero $g$.  It will be demonstrated in the later sections how the application of this approach on our MsFEM gives conveniently converging results toward the correct solution. 

\section{Application of penalization method}
\label{sec:applicationpenal}
Solving Eq. (\ref{maineq}) in $\Omega^\epsilon$ as it is often requires complex and ad-hoc grid generation methods. For highly non-periodic perforations, complicated unstructured mesh is likely what engineers would resort to. In order to confine our computations in a simple uniform Cartesian domain, we incorporate the penalization method to solve Eq. (\ref{maineq}). Henceforth, we solve instead the following 
\begin{eqnarray}
 -\nabla \cdot (\mathcal{A}^\beta \nabla u) + \vec{w}\cdot\nabla u + \sigma^\beta u = f^\beta & \textrm{in} & \Omega  \label{penal}\\
 u = g & \textrm{on} & \partial\Omega \nonumber
\end{eqnarray}
in which
\begin{eqnarray}
\mathcal{A}^\beta = \left\{ \begin{array}{rl}
 \frac{1}{h} &\mbox{ in $B^\epsilon$} \\
 \mathcal{A} &\mbox{ in $\Omega^\epsilon$}
       \end{array}  \right.,
\sigma^\beta = \left\{ \begin{array}{rl}
 \frac{1}{h^3} &\mbox{ in $B^\epsilon$} \\
 0 &\mbox{ in $\Omega^\epsilon$}
       \end{array}  \right.,
f^\beta = \left\{ \begin{array}{rl}
 0 &\mbox{ in $B^\epsilon$} \\
 f &\mbox{ in $\Omega^\epsilon$}
       \end{array}  \right..
\end{eqnarray}
Here $h$ is the width of a fine scale element used to capture highly oscillatory basis functions. We introduce the penalization coefficient $\sigma^\beta$ which forces the solution to vanish rapidly inside the perforations. Other variants of penalization methods are studied in \citep{Bruneau}.

\section{Numerical results}
\label{sec:numerical}
\subsection{Application of Bubble Functions}
In this paper, we first give numerical examples that would exhibit a case with very dense presence of perforations throughout the domain, so as to highlight the contribution of bubble function enrichment to the basis function set. In the first example, we applied bubble enrichment to a standard linearly boundary-conditioned, nodal based MsFEM without oversampling. We set aside the application of Crouzeix-Raviart in order to illustrate only the contribution of bubble functions on classical MsFEM. Taken as the computational domain is $\Omega = [0,1]^2$ with $32 \times 32$ rectangular perforations spread uniformly throughout the domain each with width of $\epsilon = 0.021875$. Dirichlet boundary conditions $u_{|\partial\Omega} = 0$ are applied in the computational domain and the source term $f = \sin (2 \pi x) \sin (2 \pi y)$ is taken. Taken as reference is the solution by standard Q1 FEM on $1024 \times 1024$ elements. 

In Figs. \ref{bubble}(a) we observe the solution of standard nodal-based MsFEM on $8 \times 8$ coarse elements with linear boundary condition without bubble function enrichment. In Figs. \ref{bubble} (b), we observe the result of the same method but with bubble function enrichment. When these two results are compared to the reference solution in Figs. \ref{bubble} (c), we notice that the one without bubble function enrichment fails to exhibit the correct solution at the interiors of coarse mesh whereas the solution with bubble function enrichment exhibits more consistency with that of the reference. In Fig. \ref{bubblebf}, we plot the standard MsFEM basis function alongside a bubble function used in this test. It clearly illustrates that with the presence of perforations this dense, the contribution of a standard MsFEM basis function inside the coarse element is insignificant. The bubble function applied in the coarse element is shown to contribute greatly to the approximation of the solution. 

In Fig. \ref{performance}, we plot the relative L2 errors with respect to the size of coarse element $H$ of the standard MsFEM with or without bubble function. We notice that the one with bubble function enrichment gives a decreasing relative error when $H$ increases away from $\epsilon$. This is of course not the behaviour exhibited by the more standard MsFEM. However, we notice that both methods shows increasing error the moment $H$ is reduced to be lower that $\epsilon$. This is due to the fact that at this region, the edges of coarse mesh start to coincide with the perforations causing incorrect solution at the interior of the MsFEM basis function. Oversampling methods  had been implemented to overcome this problem \citep{laetitia} to some degree. The contributions of Crouzeix-Raviart MsFEM as an alternative remedy to this kind of problems is reported in the next examples. 

\subsection{Application of Crouzeix-Raviart MsFEM}
\label{CRtest}
In this section, we test the Crouzeix-Raviart MsFEM with bubble functions and compare it with the standard linearly-boundary-conditioned MsFEM also with bubble functions. The test is designed to analyse the sensitivity of the methods subject to placement of perforations. The computational domain remains the same with that of the previous section. The size of each perforation is now set as $\epsilon = 0.025$. 

The methods underwent two tests: In the first test, the arrangement of the perforations is made such that none of the coarse mesh edges coincides with the them. We call this test the non-intersecting case. In the second test, the allocation of these perforations is shifted both in $x$ and $y$ direction until all coarse element edges coincide with perforations. We call this test the intersecting case. In this example, we tried to illustrate a possible worst case scenario where each and every element edges coincide with perforations at three different locations. In both cases, we implemented $8 \times 8$ coarse elements each consists of $128 \times 128$ fine elements. The reference solution is calculated using standard Q1 FEM on $1024 \times 1024$ elements. 

First, in Figs. \ref{nonintersecting}, the results of these two methods for non-intersecting case are compared with the reference solution. The results shows quantitatively good accuracies displayed by both methods. The Crouzeix-Raviart MsFEM with bubble functions records $0.11407$ L2 relative error whereas the standard MsFEM with bubble function records $0.11738$. However, in the second test, where all coarse element edges coincide with perforations, one can see in Figs. \ref{intersecting} that the standard MsFEM despite being enriched with bubble functions returns undesirable results. On the other hand, the result of Crouzeix-Raviart MsFEM with bubble functions is in good agreement with the reference recording an L2 error of $0.04269$ compared to $0.5018$ recorded by the standard MsFEM. 

To get a better understanding on why the two methods exhibit such different accuracies, we plot the basis functions of the Crouzeix-Raviart and the standard MsFEM in Figs. \ref{crbf}(a) and (b). Here one can see that the Crouzeix-Raviart basis function cope very well with perforations on the cell edges and provide natural boundary conditions around them without violating the applied constraints. By contrast, the basis function of the standard MsFEM with linear boundary condition fails to give a correct approximation in the penalized region. Again we note that several methods including oversampling methods have been introduced as remedies to this kinds of problem on standard MsFEM. Nevertheless, Crouzeix-Raviart MsFEM also has the benefit of not increasing the size of the sampled domain for constructing the MsFEM basis functions. Moreover, the exhibited natural boundary condition gives a good deal of flexibility in tackling delicate cases for it is prohibitively difficult to avoid intersections between perforations and coarse element boundaries especially when simple Cartesian mesh is implemented. In the later examples, the applicability of our method on non-periodic pattern of perforations is demonstrated. For the case of diffusion with homogeneous Dirichlet boundary condition, one can refer to the previous works in \citep{lozinskibubbleetal} where detailed comparison of performances between Crouzeix-Raviart MsFEM and other types of MsFEM including those with oversampling methods can be found. In this paper, more detailed study on the convergence behaviour of our method will be emphasized more on advection-diffusion case.

\subsection{Advection-diffusion Problems on Perforated Domain}
In this section we test our method on advection-diffusion problems on perforated domain with homogeneous boundary conditions. We implement Crouzeix-Raviart MsFEM with bubble function enrichment simply in the context of standard Galerkin approximation without any stabilizations. We reuse the computational set up done in \ref{CRtest} but with different source terms. The vector field $\vec{w}=(2y(1-x^2),-2x(1-y^2))$ set in a domain $\Omega = [-1,1]^2$. $\vec{w}$ determines a recirculating flow with streamlines $\{(x,y)|(1-x^2)(1-y^2) $ $= \textrm{constant}\}$. The source term $f(x,y)$ is defined as follows 
\begin{eqnarray}
\label{srcadv}
f(x,y) = \left\{ \begin{array}{rl}
 1 &\mbox{ if $\{-1\leq x \leq 1,0.7\leq y \leq 1\}$}\\
 1 &\mbox{ if $\{-1\leq x \leq 1,-1\leq y \leq -0.7\}$}\\
 0 &\mbox{ elsewhere. }
       \end{array}  \right.,
\end{eqnarray}
Homogeneous boundary condition $g = 0$ and diffusion parameter $\mathcal{A} = 0.03$ are applied. In Figs. \ref{convdiff} the result of both standard linearly boundary conditioned MsFEM and that of Crouzeix-Raviart MsFEM, both with bubble functions enrichments on $8 \times 8$ coarse elements, are given alongside the reference solution calculated with Q1 FEM on $1024 \times 1024$ elements. Recording a L2 relative error of $0.2287$ is the result of Crouzeix-Raviart MsFEM with bubble functions and $0.624$ recorded by the linearly boundary conditioned MsFEM with bubble functions. Clearly these results are expected given that no oversampling methods were applied. While the application of such methods is possible and could in principle improve the accuracies, this paper aims to show the possibility of using a simpler method embodied in Crouzeix-Raviart MsFEM which accommodates complex patterns of perforations without having to resort to some perforation-dependent parameters nor to enlarge the sampled domain. 

\subsection{Non-periodically perforated domain}
In this section, we test the applicability of our method on domain with non-periodic perforations. We consider two kinds of patterns of perforation as can be seen in figures \ref{perfor}. The Crouzeix-Raviart MsFEM with bubble function enrichment is implemented on $\Omega = [-1,1]^2$ domain. The first case (case a) includes 400 perforations each with width of $\epsilon = 0.025$ whereas on the second case (case b) we include 3600 perforations each with width of $\epsilon = 0.005$.  We reuse the vector field $\vec{w}=(2y(1-x^2),-2x(1-y^2))$ which determines streamlines $\{(x,y)|(1-x^2)(1-y^2) $ $= \textrm{constant}\}$. The source term (\ref{srcadv}) is also applied. On both of these cases, we utilize a diffusion coefficient of $\mathcal{A} = 0.03$. The result of the convergence tests done on these two cases is given on table \ref{deviation}. In figures (\ref{crnp1}) the contours of $u$ solved on $8 \times 8, 16 \times 16,32 \times 32,64 \times 64,$ and $128 \times 128$ elements are given alongside the reference solution on $1024 \times 1024$ solved using standard Q1-Q1 FEM. Although on $124 \times 124$ elements, the method already returns quite an identical result in comparison to the reference, the result solved on $32 \times 32$ elements is often deemed sufficient for many engineering purposes. This converging characteristic is also exhibited when solving case b, as evident from figures (\ref{crnp2}). Here the L2 relative deviations of the two cases are proportional to the values of $H/\epsilon$ as expected. 

\subsection{Application of non-homogeneous boundary condition}
Here we test the applicability of our method on solving advection-diffusion problems with non-homogeneous boundary condition.
Again we set a computational domain on $\Omega = [-1,1]^2$ where the vector field $\vec{w}=(2y(1-x^2),-2x(1-y^2))$ is set and no source term is included. Discontinuities in parts of the boundaries are introduced. At the top edge the value at the boundary is set as $u_{\partial\Omega}=1$ and $u_{\partial\Omega}=0$ everywhere else.  Randomly placed $100$ perforations are considered each with width of $\epsilon = 0.04$ as shown in figure (\ref{perfor2}). In figures (\ref{crnp3}) the contours of $u$ solved on $8 \times 8, 16 \times 16,32 \times 32,64 \times 64,$ and $128 \times 128$ elements are given alongside the reference solution on $1024 \times 1024$ solved using standard Q1-Q1 FEM. In table \ref{deviation_nonbc}, it is shown that the method returns grid converging results toward the reference solution as exhibited in previous tests with homogeneous boundary conditions.

%%%%%%%%%%%%%%%%%%%%%%%%%%%%%%%%%%TABLES%%%%%%%%%%%%%%%%%%%%%%%%%%%%%%%%%%%%%%%%%
\begin{table}[hbt]
\begin{center}
\begin{tabular}{{l}{l}{l}{l}{l}{c}}
\hline\hline
 && \multicolumn{2}{c}{$H/\epsilon$}& \multicolumn{2}{c}{L2} \\ 
Config.      			& H         &   case a   &	case b	& case a &	case b\\
\hline
$8\times 8$     	& 0.25 		&10		&50 		&0.273 & 0.346\\
$16\times 16$     	& 0.125 	&5		&25 		&0.265 & 0.337\\
$32\times 32$     	& 0.0625 	&2.5	&12.5 		&0.140 & 0.321\\
$64\times 64$   	& 0.03125 	&1.25	&6.25 		&0.098 & 0.284\\
$128\times 128$   	& 0.015625 	&0.625	&3.125 		&0.031 & 0.148\\
\hline\hline
\end{tabular}
\end{center}
\caption{Deviation from reference solution for case a ($\epsilon = 0.025$) and case b ($\epsilon = 0.005$)}
\label{deviation}
\end{table}
%%%%%%%%%%%%%%%%%%%%%%%%%%%%%%%%%%%%%%%%%%%%%%%%%%%%%%%%%%%%%%%%%%%%%%%%%%%%%%%%%%
%%%%%%%%%%%%%%%%%%%%%%%%%%%%%%%%%%%%%%%%%%%%%%%%%%%%%%%%%%%%%%%%%%%%%%%%%%%%%%%%%%
\begin{table}[hbt]
\begin{center}
\begin{tabular}{{l}{l}{l}{l}{l}{c}}
\hline\hline
Config.      			& H         &   {$H/\epsilon$}	& {L2}\\
\hline
$8\times 8$     	& 0.25 		&6.25		 		&0.487 \\
$16\times 16$     	& 0.125 	&3.125		 		&0.206 \\
$32\times 32$     	& 0.0625 	&1.5625	 			&0.073 \\
$64\times 64$   	& 0.03125 	&0.78125	 		&0.027 \\
$128\times 128$   	& 0.015625 	&0.390625	 		&0.013 \\
\hline\hline
\end{tabular}
\end{center}
\caption{Deviation from reference solution (with non-homogeneous B.C. and $\epsilon = 0.04$)}
\label{deviation_nonbc}
\end{table}
%%%%%%%%%%%%%%%%%%%%%%%%%%%%%%%%%%%%%%%%%%%%%%%%%%%%%%%%%%%%%%%%%%%%%%%%%%%%%%%%%%

\section{Concluding remarks}
\label{sec:concluding}
In this paper, the feasibility of Crouzeix-Raviart MsFEM with bubble function enrichments for solving diffusion and advection-diffusion problems in perforated media through means of penalization methods have been demonstrated without much major constraints. The resulting method allows us to address multiscale problems with inconvenient patterns of perforations and still obtain accurate solutions between perforations. Although in the given examples, the diffusion coefficient $\mathcal{A}$ are taken as constants, Crouzeix-Raviart MsFEM has been shown to be able to solve highly oscillatory problems \citep{clbetal}. Crouzeix-Raviart MsFEM with bubble function enrichment has shown good performance in comparison with more conventional MsFEMs especially in as far as insensitivity to size and placements of perforations is concerned. We also include the cases for non-periodic perforations where the robustness of our method is tested in more realistic circumstances. 
%%%% Acknowledgments %%%%%%%%

%%%% Bibliography  %%%%%%%%%%
\bibliographystyle{model6-num-names}
\bibliography{crmsfem}

\end{document}